\documentclass[10pt]{article}
\usepackage{amsmath}
\usepackage{amsfonts}
\usepackage{amssymb}
\usepackage{graphicx}
\usepackage{mathrsfs}
\usepackage{xcolor}
\usepackage{verbatim}
\usepackage{mathrsfs}
\usepackage[body={15.5cm,21cm}, top=3cm]{geometry}
\DeclareMathOperator*{\esssup}{ess\,sup}
\usepackage{paralist}

\allowdisplaybreaks[4]
\providecommand{\U}[1]{\protect \rule{.1in}{.1in}}
\newtheorem{theorem}{Theorem}[section]

\newenvironment{proof}[1][Proof]{\noindent \textbf{#1.} }{\  \rule{0.5em}{0.5em}}
\allowdisplaybreaks[2]
\begin{document}
\title{Existence and uniqueness of solutions for multi-dimensional reflected BSDEs with diagonally quadratic generators}
\author{Yuyang Chen\thanks{School of Mathematical Sciences, Shanghai Jiao Tong University, China (cyy0032@sjtu.edu.cn)}
\and
Peng Luo \thanks{School of Mathematical Sciences, Shanghai Jiao Tong University, China (peng.luo@sjtu.edu.cn). Financial support from the National Natural Science Foundation of China (Grant No. 12101400) is gratefully acknowledged.}}

\maketitle
\begin{abstract}
In this paper, we study multi-dimensional reflected backward stochastic differential equations with diagonally quadratic generators. Using the comparison theorem for diagonally quadratic BSDEs which is established recently in \cite{19}, we obtain the existence and uniqueness of a solution by a penalization method. Moreover, we provide a comparison theorem.
\end{abstract}


\textbf{MSC-classification}: 60H10, 60H30.

\section{Introduction}
Reflected backward stochastic differential equations (RBSDEs for short) was first introduced by El Karoui et al. \cite{1}. They obtained existence and uniqueness of solutions for Lipschitz generators and established its connection with obstacle problem of related PDEs. Matoussi \cite{3} proved the existence of maximal and minimal solutions for RBSDEs with continuous generators having a linear growth. Lepeltier et al. \cite{4} obtained the existence and uniqueness of the solution by using an approximation procedure when the generator is Lipschitz in $z$ and satisfies a monotone condition in $y$. $L^{p}$ solutions of RBSDEs were studied by Klimsiak \cite{5}. All these results are one-dimesional. We refer to Gegout-Petit and Pardoux \cite{GP}, Chassagneux et al. \cite{14}, Wu and Xiao \cite{18} for results in multi-dimensional case. 

The case where the generator is allowed to have quadratic growth in $z$ was considered in Kobylanski et al. \cite{11}. They obtained the existence of a maximal and minimal bounded solution for RBSDE when the generator has a superlinear growth in $y$ and quadratic growth in $z$. This result was further extended by Lepeltier and Xu \cite{12} for unbounded terminal conditions. Xu \cite{13} proved the existence and uniqueness of a solution for RBSDE when the generator satifies a monotone condition in $y$ and non-Lipschitz condition in $z$.

Motivated by recent developments of solvability of multi-dimesional quadratic BSDEs (see \cite{Te,17,JKL,Luo,CN, XZ}), we consider multi-dimensional RBSDEs with diagonally quadratic generators. By a standard technique, we introcude the penalized BSDEs and provide some delicated a priori estimates for their solutions. Using a comparison result for diagonally quadratic BSDEs obtained in \cite{19}, we first show the monotone convergence of the sequence of value processes of the penalized BSDEs. Relyling on properties of BMO martingales, we further show that solutions of the penalized BSDEs converge in some suitable spaces, from which we obtain the existence of a solution for the multi-dimeional RBSDE with a diagonally quadratic generator. Moreover, we verify that the solution belongs to some finer space which allows us to show the uniqueness by a standard argument. On the other hand, we provide a representation result for the value process in terms of an optimal stopping problem. Finally, we give a comparison theorem for multi-dimensional RBSDEs with diagonally quadratic generators. The paper is organized as follows. In the next section, we state our setup and main results. The existence and uniqueness results are given in Theorem 2.1. Theorem 2.2 provides a representation of the value process by a optimal stopping problem, while Theorem 2.3 gives the comparison result.

\section{Main results}
Consider the following reflected backward stochastic differential equations(RBSDEs) with one continuous barrier:
\begin{equation}
\begin{aligned}
&Y_{t}^{i}=\xi^{i}+\int_{t}^{T}f^{i}(s,Y_{s},Z_{s}^{i})ds+K_{T}^{i}-K_{t}^{i}-\int_{t}^{T}Z_{s}^{i}dW_{s},\\
&Y_{t}^{i}\geq S_{t}^{i},~\int_{0}^{T}(Y_{t}^{i}-S_{t}^{i})dK_{t}^{i}=0
\end{aligned}
\end{equation}
for $i=1,\cdots,n,~t\in[0,T]$, where $K$ is a $n$-dimensional continuous and nondecreasing process with $K_{0}=0$. Here, $W=(W_{t})_{t\geq0}$ is a $d$-dimensional standard Brownian motion defined on some probability space $(\Omega,\mathcal{F},P)$. Denote by $\{\mathcal{F}_{t},0\leq t\leq T\}$ the augmented natural filtration of W. Equalities and inequalities between random variables and processes are understood in the $P$-a.s. and $P\otimes dt$-a.e. sense, respectively. The Euclidean norm is denoted by $|\cdot|$ and for a random variable $\eta$, $\|\eta\|_{\infty}$ denotes the $L^{\infty}$-norm of $\eta$, i.e., $\|\eta\|_{\infty}:=\mathop{\esssup}\limits_{\omega}|\eta(\omega)|$. For $x,y\in\mathbb{R}^{n},~x\leq y$ is understood component-wisely, i.e., $x\leq y$ if and only if $x^{i}\leq y^{i}$ for all $i=1,\cdots,n$. For $p\geq2$, we denote by
\begin{itemize}
  \item $\mathcal{H}^{p}(\mathbb{R}^{n})$\ the space of $n$-dimensional predictable processes $Y$ on [0,T] such that\\
  \centerline{$\|Y\|_{\mathcal{H}^{p}(\mathbb{R}^{n})}:=E\left[\left(\int_{0}^{T}|Y_{t}|^{2}dt\right)^{\frac{p}{2}}\right]^{\frac{1}{p}}<\infty$;}\\
  \item $\mathcal{S}^{p}(\mathbb{R}^{n})$\ the space of $n$-dimensional predictable processes $Y$ on [0,T] such that\\
  \centerline{$\|Y\|_{\mathcal{S}^{p}(\mathbb{R}^{n})}:=E\left[\left(\sup\limits_{0\leq t\leq T}|Y_{t}|^{2}\right)^{\frac{p}{2}}\right]^{\frac{1}{p}}<\infty$;}\\
  \item $\mathcal{S}^{\infty}(\mathbb{R}^{n})$\ the space of $n$-dimensional predictable processes $Y$ on [0,T] such that\\
  \centerline{$\|Y\|_{\mathcal{S}^{\infty}(\mathbb{R}^{n})}:=\left\|\sup\limits_{0\leq t\leq T}|Y_{t}|\right\|_{\infty}<\infty$;}\\
  \item $L^{\infty}(\mathbb{R}^{n})$\ the space of $n$-dimensional $\mathcal{F}_{T}$-measurable random variables $\xi$ such that\\
  \centerline{$\|\xi\|_{\infty}<\infty$.}
\end{itemize}
Let $\mathcal{T}$ be the set of all stopping times with values in $[0,T]$. For any uniformly integrable martingale $M$ with $M_{0}=0$, we set
\begin{align*}
\|M\|_{BMO_{2}}:=\sup_{\tau\in\mathcal{T}}\|E[|M_{T}-M_{\mathcal{\tau}}|^{2}|\mathcal{F}_{\mathcal{T}}]^{\frac{1}{2}}\|_{\infty}.
\end{align*}
The class $\{M:\|M\|_{BMO_{2}}<\infty\}$ is denoted by $BMO$. For $(\alpha\cdot W)_{t}:=\int_{0}^{t}\alpha_{s}dW_{s}$ in $BMO$, the corresponding stochastic exponential is denoted by $\mathcal{E}_{t}(\alpha\cdot W)$.\\
In this paper, we make the following assumptions. Let $C$ be a positive constant.
\begin{flushleft}
$(\mathscr{A}1)$For $i=1,\cdots,n$, the function $f^{i}:\Omega\times[0,T]\times\mathbb{R}^{n}\times\mathbb{R}^{d}\rightarrow\mathbb{R}$ satisfies that $f^{i}(\cdot,y,z^{i})$ is adapted for each $y\in\mathbb{R}^{n}$ and $z^{i}\in\mathbb{R}^{d}$. It holds that
\begin{align*}
&|f^{i}(t,y,z^{i})|\leq C(1+|y|+|z^{i}|^{2}),\\
&|f^{i}(t,y,z^{i})-f^{i}(t,\overline{y},\overline{z}^{i})|\leq C|y-\overline{y}|+C(1+|z^{i}|+|\overline{z}^{i}|)|z^{i}-\overline{z}^{i}|
\end{align*}
for $y,\overline{y}\in\mathbb{R}^{n}$, $z^{i},\overline{z}^{i}\in\mathbb{R}^{d}$.
\end{flushleft}
\begin{flushleft}
$(\mathscr{A}2)$For $t\in[0,T],~i=1,\cdots,n$, it holds that
\begin{align*}
f^{i}(t,y,z^{i})\leq f^{i}(t,\overline{y},z^{i})
\end{align*}
for any $y,\overline{y}$ satisfying $y^{i}=\overline{y}^{i},y^{j}\leq \overline{y}^{j},j\neq i$.
\end{flushleft}
\begin{flushleft}
$(\mathscr{A}3)$$S$ is a continuous progressively measurable real-value process satisfying
\begin{align*}
\|S^{+}\|_{\mathcal{S}^{\infty}(\mathbb{R}^{n})}\leq C.
\end{align*}
\end{flushleft}
\begin{flushleft}
$(\mathscr{A}4)$The terminal condition $\xi$ satisfies
\begin{align*}
\|\xi\|_{L^{\infty}(\mathbb{R}^{n})}\leq C.
\end{align*}
\end{flushleft}

\begin{theorem}
Let $(\mathscr{A}1)-(\mathscr{A}4)$ be satisfied, then RBSDE (1) has a unique solution $(Y,Z,K)$ such that $(Y,Z\cdot W,K)\in\mathcal{S}^{\infty}(\mathbb{R}^{n})\times BMO\times\mathcal{S}^{p}(\mathbb{R}^{n})$ for any $p\geq2$.
\end{theorem}
\begin{proof}
We divide the proof into four steps:\\
$\mathbf{Step\ 1}$: For any $t\in[0,T],~k\geq0,~i=1,\cdots,n$, we consider the following RBSDEs:
\begin{align}
Y_{t}^{i,(k)}=\xi^{i}+\int_{t}^{T}f^{i}(s,Y_{s}^{(k)},Z_{s}^{i,(k)})ds+k\int_{t}^{T}(Y_{s}^{i,(k)}-S_{s}^{i})^{-}ds-\int_{t}^{T}Z_{s}^{i,(k)}dW_{s}.
\end{align}
For fixed $k$, RBSDE (2) has a unique solution $(Y^{(k)},Z^{(k)})$ such that $(Y^{(k)},Z^{(k)}\cdot W)\in\mathcal{S}^{\infty}(\mathbb{R}^{n})\times BMO$ due to \cite[Theorem 2.3]{17}.\\
Define $K_{t}^{i,(k)}=k\int_{0}^{t}(Y_{s}^{i,(k)}-S_{s}^{i})^{-}ds$ and let
\begin{align*}
f^{i}(t,Y_{t}^{(k)},Z_{t}^{i,(k)})-f^{i}(t,0,0)=\alpha^{i,k}_{t} Y_{t}^{(k)}+\beta^{i,k}_{t} Z_{t}^{i,(k)},~|\alpha^{i,k}_{t}|\leq C,~|\beta^{i,k}_{t}|\leq C(1+|Z_{t}^{i,(k)}|).
\end{align*}
Then $W_{t}^{i,k}:=W_{t}-\int_{0}^{t}\beta_{s}^{i,k}ds$ is a Brownian motion under the equivalent probability measure $P^{i,k}$ defined by
\begin{align*}
dP^{i,n}:=\mathcal{E}(\beta^{i,k}\cdot W)_{0}^{T}dP\\
\end{align*}
and
\begin{align*}
Y_{t}^{i,(k)}=\xi^{i}+\int_{t}^{T}f^{i}(s,0,0)ds+\int_{t}^{T}\alpha^{i,k}_{s} Y_{s}^{(k)}ds+K_{T}^{i,(k)}-K_{t}^{i,(k)}-\int_{t}^{T}Z_{s}^{i,(k)}dW^{i,k}_{s}.
\end{align*}
Use It\^o's formula to $|Y_{t}^{i,(k)}|^{2}$ under the conditional expectation $E_{t}^{i,k}$ and we will get
\begin{equation*}
\begin{aligned}
&|Y_{t}^{i,(k)}|^{2}+E_{t}^{i,k}\int_{t}^{T}|Z_{s}^{i,(k)}|^{2}ds\\
&=E_{t}^{i,k}|\xi^{i}|^{2}+2E_{t}^{i,k}\int_{t}^{T}f^{i}(s,0,0)Y_{s}^{i,(k)}ds+2E_{t}^{i,k}\int_{t}^{T}(\alpha^{i,k}_{s} Y_{s}^{(n)})Y_{s}^{i,(k)}ds+2E_{t}^{i,k}\int_{t}^{T}Y_{s}^{i,(k)}dK_{s}^{i,(k)}.\\
\end{aligned}
\end{equation*}
We calculate
\begin{align*}
&E_{t}^{i,k}\left[\left(K_{T}^{i,(k)}-K_{t}^{i,(k)}\right)^{2}\right]\\
&\leq5\{|Y_{t}^{i,(k)}|^{2}+E_{t}^{i,k}|\xi^{i}|^{2}+E_{t}^{i,k}|\int_{t}^{T}f^{i}(s,0,0)ds|^{2}+E_{t}^{i,k}|\int_{t}^{T}\alpha^{i,k}_{s}Y_{s}^{(n)}ds|^{2}+E_{t}^{i,k}\int_{t}^{T}|Z^{i,(k)}_{s}|^{2}ds\}\\
&\leq5|Y_{t}^{i,(k)}|^{2}+5C^{2}+5C^{2}(T-t)^{2}+5C^{2}(T-t)^{2}\sum_{j=1}^{l}E_{t}^{i,k}\left[\sup\limits_{t\leq s\leq T}|Y_{s}^{j,(n)}|^{2}\right]+5E_{t}^{i,k}\int_{t}^{T}|Z^{i,(k)}_{s}|^{2}ds
\end{align*}
and then
\begin{align*}
&|Y_{t}^{i,(k)}|^{2}+E_{t}^{i,k}\int_{t}^{T}|Z_{s}^{i,(k)}|^{2}ds\\
&\leq E_{t}^{i,k}|\xi^{i}|^{2}+E_{t}^{i,k}\int_{t}^{T}|f^{i}(s,0,0)|^{2}ds+E_{t}^{i,k}\int_{t}^{T}|\alpha^{i,k}_{s}Y_{s}^{(k)}|^{2}ds+2E_{t}^{i,k}\int_{t}^{T}|Y_{s}^{i,(k)}|^{2}ds\\
&\quad+2E_{t}^{i,k}\int_{t}^{T}S_{s}^{i,(k)}dK_{s}^{i,(k)}\\
&\leq C^{2}+C^{2}(T-t)+C^{2}(T-t)\sum_{j=1}^{n}E_{t}^{i,k}\left[\sup\limits_{t\leq s\leq T}|Y_{s}^{j,(k)}|^{2}\right]+2(T-t)E_{t}^{i,k}\left[\sup\limits_{t\leq s\leq T}|Y_{s}^{i,(k)}|^{2}\right]\\
&\quad+\frac{1}{\alpha}E_{t}^{i,k}|\sup\limits_{t\leq s\leq T}(S_{s}^{i})^{+}|^{2}+\alpha E_{t}^{i,k}\left[(K_{T}^{i,(k)}-K_{t}^{i,(k)})^{2}\right]\\
&\leq[1+(T-t)+5\alpha+5\alpha(T-t)^{2}+\frac{1}{\alpha}]C^{2}+5\alpha|Y_{t}^{i,(k)}|^{2}+5\alpha E_{t}^{i,k}\int_{t}^{T}|Z^{i,(k)}_{s}|^{2}ds\\
&\quad+\left[5\alpha C^{2}(T-t)+C^{2}+2\right](T-t)\sum_{j=1}^{n}E_{t}^{i,k}\left[\sup\limits_{t\leq s\leq T}|Y_{s}^{j,(n)}|^{2}\right].
\end{align*}
Let $\alpha=\frac{1}{10}$, we have
\begin{align*}
&\frac{1}{2}\|\sup_{t\leq s\leq T}|Y_{s}^{i,(k)}|^{2}\|_{\infty}\\
&\leq[\frac{23}{2}+\frac{1}{2}(T-t)^{2}+(T-t)]C^{2}+[\frac{1}{2}C^{2}(T-t)+C^{2}+2](T-t)\sum_{j=1}^{n}\|\sup_{t\leq s\leq T}|Y_{s}^{j,(k)}|^{2}\|_{\infty},
\end{align*}
then sum up the equation above respect to $i$ from $1$ to $n$, it holds that
\begin{align*}
&\frac{1}{2}\sum_{i=1}^{n}\|\sup_{t\leq s\leq T}|Y_{s}^{i,(k)}|^{2}\|_{\infty}\\
&\leq[\frac{23}{2}+\frac{1}{2}(T-t)^{2}+(T-t)]nC^{2}+[\frac{1}{2}C^{2}(T-t)+C^{2}+2](T-t)n\sum_{i=1}^{n}\|\sup_{t\leq s\leq T}|Y_{s}^{i,(k)}|^{2}\|_{\infty}.
\end{align*}
Let
\begin{equation*}
\begin{aligned}
\delta=\frac{1}{2}-[\frac{1}{2}C^{2}(T-t)+C^{2}+2](T-t)n,
\end{aligned}
\end{equation*}
when $T-t=\varepsilon>0$ is small enough, we have $0<\delta<\frac{1}{2}$ and
\begin{equation*}
\begin{aligned}
\|\sup\limits_{t\leq s\leq T}|Y_{s}^{i,(k)}|^{2}\|_{\infty}\leq[\frac{23}{2}+\frac{1}{2}(T-t)^{2}+(T-t)]\frac{nC^{2}}{\delta}.
\end{aligned}
\end{equation*}
Let
\begin{equation*}
\begin{aligned}
C'=[\frac{23}{2}+\frac{1}{2}(T-t)^{2}+(T-t)]^{\frac{1}{2}}(\frac{nC^{2}}{\delta})^{\frac{1}{2}}
\end{aligned}
\end{equation*}
we can get
\begin{equation*}
\begin{aligned}
\|\sup\limits_{T-2\varepsilon\leq s\leq T-\varepsilon}|Y_{s}^{i,(k)}|^{2}\|_{\infty}\leq[\frac{23}{2}+\frac{1}{2}(T-t)^{2}+(T-t)]\frac{nC'^{2}}{\delta}
\end{aligned}
\end{equation*}
with the same method above. Since the terminal time $T$ is finite, after a finite step iteration, e.g. $m$ steps, we have
\begin{equation*}
\begin{aligned}
\|\sup\limits_{0\leq s\leq T}|Y_{s}^{i,(k)}|^{2}\|_{\infty}\leq[\frac{23}{2}+\frac{1}{2}(T-t)^{2}+(T-t)]^{m}(\frac{n}{\delta})^{m}C^{2}.
\end{aligned}
\end{equation*}
For convenience, we assume $\|Y^{(k)}\|_{\mathcal{S}^{\infty}(\mathbb{R}^{n})}\leq C^{*}$ where $C^{*}$ is independent of $k$.\\
Let $\phi(x)=e^{-4Cx}$ and we can calculate $\phi'(x)=-4Ce^{-4Cx},~\phi''(x)=16C^{2}e^{-4Cx}$. Use It\^o's formula to $\phi(Y_{t}^{i,(k)})$ under the conditional expectation $E_{t}$,
\begin{align*}
&\phi(Y_{t}^{i,(k)})+\frac{1}{2}E_{t}\int_{t}^{T}\phi''(Y_{s}^{i,(k)})|Z_{s}^{i,(n)}|^{2}ds\\
&=E_{t}\phi(\xi^{i})+E_{t}\int_{t}^{T}\phi'(Y_{s}^{i,(k)})f^{i}(s,Y_{s}^{(k)}Z_{s}^{i,(k)})ds+E_{t}\int_{t}^{T}\phi'(Y_{s}^{i,(k)})dK_{s}^{i,(k)}.
\end{align*}
Since $\phi'(Y_{t}^{i,(k)})\leq0,~K_{t}^{i,(k)}\geq0,~|f^{i}(t,Y_{t}^{(k)}Z_{t}^{i,(k)})|\leq C(1+|Y_{t}^{(k)}|+|Z_{t}^{i,(k)}|^{2})$, we have
\begin{equation*}
\begin{aligned}
\phi(Y_{t}^{i,(k)})+\frac{1}{2}E_{t}\int_{t}^{T}\phi''(Y_{s}^{i,(k)})|Z_{s}^{i,(k)}|^{2}ds&\leq E_{t}\phi(\xi^{i})+E_{t}\int_{t}^{T}C|\phi'(Y_{s}^{i,(k)})|(1+|Y_{t}^{(k)}|+|Z_{t}^{i,(k)})|^{2})ds\\
\end{aligned}
\end{equation*}
and furthermore
\begin{equation}
\begin{aligned}
4C^{2}e^{-4CC^{*}}E_{t}\int_{t}^{T}|Z_{s}^{i,(k)}|^{2}ds
&\leq e^{4C^{2}}+4C^{2}(1+C^{*})Te^{4CC^{*}}
\end{aligned}
\end{equation}
which implies there exists $A>0$ independent of $k$, such that $\|Z^{(k)}\cdot W\|_{BMO}\leq A$.\\
Then it is easily seen that for any $p\geq2$,
\begin{align*}
&E\left[\left(K_{T}^{i,(k)}\right)^{p}\right]\\
&\leq C_{p}\left\{E|Y_{0}^{i,(k)}|^{p}+E|\xi^{i}|^{p}+E|\int_{0}^{T}f^{i}(s,0,0)ds|^{p}\right.\\
&\quad+\left.E|\int_{0}^{T}\alpha_{s}^{i,k}Y_{s}^{(n)}ds|^{p}+E|\int_{0}^{T}\beta_{s}^{i,k}Z_{s}^{i,(k)}ds|^{p}+E|\int_{0}^{T}Z_{s}^{i,(k)}dW_{s}|^{p}\right\}\\
&\leq C_{p}\left\{(C^{*})^{p}+C^{p}+(CT)^{p}+(CC^{*}T)^{p}\right.\\
&\quad+\left.K^{p}E\left[\left(\int_{0}^{T}(1+|Z_{s}^{i,(k)}|)^{2}ds\right)^{p}\right]^{\frac{1}{2}}E\left[\left(\int_{0}^{T}|Z_{s}^{i,(k)}|^{2}ds\right)^{p}\right]^{\frac{1}{2}}+E\left[\left(\int_{0}^{T}|Z_{s}^{i,(k)}|^{2}ds\right)^{\frac{p}{2}}\right]\right\},
\end{align*}
which implies there exists $B_{p}>0$ independent of $k$, such that $E\left[\left(K_{t}^{(k)}\right)^{p}\right]^{\frac{1}{p}}\leq E\left[\left(K_{T}^{(k)}\right)^{p}\right]^{\frac{1}{p}}\leq B_{p}$ for any $p\geq2,~t\in[0,T]$.\\
\\
$\mathbf{Step\ 2}$: We prove $\lim\limits_{k\rightarrow\infty}E[\sup\limits_{0\leq t\leq T}|(Y_{t}^{i,(k)}-S_{t}^{i})^{-}|^{p}]=0$ for any $p\geq2$.\\
Inspired by \cite[Lemma 4.2]{21}, set
\begin{align*}
y_{t}^{i,(k)}=\xi^{i}+\int_{t}^{T}f^{i}(s,0,0)ds+\int_{t}^{T}\alpha_{s}^{i,k}Y_{s}^{(k)}ds+k\int_{t}^{T}(S_{s}^{i}-y_{s}^{i,(k)})ds-\int_{t}^{T}z_{s}^{i,(k)}dW_{s}^{i,k}
\end{align*}
and we immediately obtain $y_{t}\leq Y_{t},~t\in[0,T]$ from \cite[Theorem 2.2]{19}.\\
Let $\nu\in[0,T]$ be a stopping time, then
\begin{align*}
y_{\nu}^{i,(k)}=E_{\nu}^{i,k}[e^{-k(T-\nu)}\xi^{i}+\int_{\nu}^{T}e^{-ks}(f^{i}(s,0,0)+\alpha_{s}^{i,k}Y_{s}^{(k)})ds+\int_{\nu}^{T}ke^{-k(s-\nu)}S_{s}^{i}ds].
\end{align*}
It is easily seen that
\begin{align*}
E_{\nu}|(e^{-k(T-\nu)}\xi^{i}+\int_{\nu}^{T}ke^{-k(s-\nu)}S_{s}^{i}ds)-(\xi^{i}1_{\{\nu=T\}}+S_{\nu}^{i}1_{\{\nu<T\}})|^{p}\rightarrow0\ as\ k\rightarrow\infty.
\end{align*}
We note that there exists $p(i,k)>1$, such that
\begin{equation*}
\begin{aligned}
E_{\nu}\left(\frac{dP^{i,k}}{dP}\right)^{p(i,k)}<\infty
\end{aligned}
\end{equation*}
from \cite[Theorem 3.1]{20}, where $E_{\nu}$ is a conditional expectation for an arbitrary stopping time $\nu$. Let $q(i,k)>1$, such that
\begin{equation*}
\begin{aligned}
\frac{1}{p(i,k)}+\frac{1}{q(i,k)}=1.
\end{aligned}
\end{equation*}
Then for any $p\geq2$,
\begin{align*}
&E_{\nu}^{i,k}|(e^{-k(T-\nu)}\xi^{i}+\int_{\nu}^{T}ke^{-k(s-\nu)}S_{s}^{i}ds)-(\xi^{i}1_{\{\nu=T\}}+S_{\nu}^{i}1_{\{\nu<T\}})|^{p}\\
&=E_{\nu}\left(\frac{dP^{i,k}}{dP}\right)|(e^{-k(T-\nu)}\xi^{i}+\int_{\nu}^{T}ke^{-k(s-\nu)}S_{s}^{i}ds)-(\xi^{i}1_{\{\nu=T\}}+S_{\nu}^{i}1_{\{\nu<T\}})|^{p}\\
&\leq E_{\nu}\left[\left(\frac{dP^{i,k}}{dP}\right)^{p(i,k)}\right]^{\frac{1}{p(i,k)}}E_{\nu}\left[|(e^{-k(T-\nu)}\xi^{i}+\int_{\nu}^{T}ke^{-k(s-\nu)}S_{s}^{i}ds)-(\xi^{i}1_{\{\nu=T\}}+S_{\nu}^{i}1_{\{\nu<T\}})|^{pq(i,k)}\right]^{\frac{1}{q(i,k)}}.
\end{align*}
Similarly, for any $p\geq2$, we can calculate
\begin{align*}
&E_{\nu}^{i,k}|\int_{\nu}^{T}e^{-ks}(f^{i}(s,0,0)+\alpha_{s}^{i,k}Y_{s}^{(k)})ds|^{p}\\
&=E_{\nu}\left(\frac{dP^{i,k}}{dP}\right)|\int_{\nu}^{T}e^{-ks}(f^{i}(s,0,0)+\alpha_{s}^{i,k}Y_{s}^{(k)})ds|^{p}\\
&\leq E_{\nu}\left[\left(\frac{dP^{i,k}}{dP}\right)^{p(i,k)}\right]^{\frac{1}{p(i,k)}}E_{\nu}\left[|\int_{\nu}^{T}e^{-ks}(f^{i}(s,0,0)+\alpha_{s}^{i,k}Y_{s}^{(k)})ds|^{pq(i,k)}\right]^{\frac{1}{q(i,k)}}\\
&\leq E_{\nu}\left[\left(\frac{dP^{i,k}}{dP}\right)^{p(i,k)}\right]^{\frac{1}{p(i,k)}}E_{\nu}\left[|\int_{\nu}^{T}e^{-2ks}ds\int_{\nu}^{T}(f^{i}(s,0,0)+\alpha_{s}^{i,k}Y_{s}^{(k)})^{2}ds|^{\frac{pq(i,k)}{2}}\right]^{\frac{1}{q(i,k)}}\\
&\leq\left(\frac{1}{2k}\right)^{\frac{p}{2}}E_{\nu}\left[\left(\frac{dP^{i,k}}{dP}\right)^{p(i,k)}\right]^{\frac{1}{p(i,k)}}E_{\nu}\left[|\int_{\nu}^{T}(C+CC^{*})^{2}ds|^{\frac{pq(i,k)}{2}}\right]^{\frac{1}{q(i,k)}}\\
&\leq\left(\frac{TC^{2}(1+C^{*})^{2}}{2k}\right)^{\frac{p}{2}}E_{\nu}\left[\left(\frac{dP^{i,k}}{dP}\right)^{p(i,k)}\right]^{\frac{1}{p(i,k)}}
\end{align*}
Consequently, with the arbitrariness of $p$, it is easily seen that
\begin{align*}
\lim\limits_{k\rightarrow\infty}E|y_{\nu}^{i,(k)}-\xi^{i}1_{\{\nu=T\}}-S_{\nu}^{i}1_{\{\nu<T\}}|^{p}=0,~\forall p\geq2
\end{align*}
and then
\begin{align*}
Y_{\nu}^{i}=\lim\limits_{k\rightarrow\infty}Y_{\nu}^{i,(k)}\geq \lim\limits_{k\rightarrow\infty}y_{\nu}^{i,(k)}=S_{\nu}^{i}.
\end{align*}
By the optional section theorem, we have
\begin{align*}
Y_{t}^{i}\geq S_{t}^{i},~t\in[0,T].
\end{align*}
We can deduce
\begin{align*}
\lim\limits_{k\rightarrow\infty}\sup\limits_{0\leq t\leq T}(Y_{t}^{i,(k)}-S_{t}^{i})^{-}=0
\end{align*}
from $\lim\limits_{k\rightarrow\infty}(Y_{t}^{i,(k)}-S_{t}^{i})^{-}=0,~t\in[0,T]$ and by Dini's theorem.\\
Due to $(Y_{t}^{i,(k)}-S_{t}^{i})^{-}\leq |Y_{t}^{i,(k)}|+(S_{t}^{i})^{+}\leq C^{*}+C$ and the dominated convergence theorem, we have
\begin{align*}
\lim\limits_{k\rightarrow\infty}E[\sup\limits_{0\leq t\leq T}|(Y_{t}^{i,(k)}-S_{t}^{i})^{-}|^{p}]=0,~\forall p\geq2.
\end{align*}
\\
$\mathbf{Step\ 3}$: We prove that $(Y^{(k)},Z^{(k)},K^{(k)})$ converges to $(Y,Z,K)$ in $\mathcal{S}^{p}(\mathbb{R}^{n})\times\mathcal{H}^{p}(\mathbb{R}^{n})\times\mathcal{S}^{p}(\mathbb{R}^{n})$ as $k\rightarrow\infty$ for any $p\geq2$.\\
We define $f_{k}^{i}(t,y,z^{i})=f^{i}(t,y,z^{i})+k(y^{i}-S_{t}^{i})^{-},~k\geq0$. We can deduce
\begin{align*}
Y_{t}^{i,(k)}\leq Y_{t}^{i,(k+1)}\leq C^{*},~t\in[0,T]
\end{align*}
from $f_{k}^{i}(t,y,z^{i})\leq f_{k+1}^{i}(t,y,z^{i}),~t\in[0,T]$ and \cite[Theorem 2.2]{19}. Therefore, there exists an adapted process $Y$ such that $\|Y^{(k)}-Y\|_{\mathcal{H}^{p}(\mathbb{R}^{n})}\rightarrow0$ as $k\rightarrow\infty$ for any $p\geq2$ by the monotonic convergence theorem. What's more, $\lim\limits_{k,h\rightarrow\infty}E\int_{0}^{T}|Y_{t}^{(k)}-Y_{t}^{(h)}|^{p}dt=0$, for any $p\geq2$.\\
For $k,h\geq0,~i=1,\cdots,n$, denote
\begin{align*}
&\Delta Y^{k,h}=Y^{(k)}-Y^{(h)},\Delta Y^{i,k,h}=Y^{i,(k)}-Y^{i,(h)},\\
&\Delta Z^{i,k,h}=Z^{i,(k)}-Z^{i,(h)},\Delta K^{i,k,h}=K^{i,(k)}-K^{i,(h)},\\
&f^{i}(t,Y_{t}^{(k)},Z_{t}^{i,(k)})-f^{i}(t,Y_{t}^{(h)},Z_{t}^{i,(h)})=\alpha^{i,k,h}_{t} \Delta Y_{t}^{k,h}+\beta^{i,k,h}_{t} \Delta Z_{t}^{i,k,h},\\
&|\alpha^{i,k,h}_{t}|\leq C,~\ |\beta^{i,k,h}_{t}|\leq C(1+|Z_{t}^{i,(k)}|+|Z_{t}^{i,(k)}|).
\end{align*}
Then $W_{t}^{i,k,h}:=W_{t}-\int_{0}^{t}\beta_{s}^{i,k,h}ds$ is a Brownian motion under the equivalent probability measure $P^{i,k,h}$ defined by
\begin{align*}
dP^{i,k,h}:=\mathcal{E}(\beta^{i,k,h}\cdot W)_{0}^{T}dP.
\end{align*}
Use It\^o's formula to $|\Delta Y_{t}^{i,k,h}|^{p}$ for any $p\geq2$ under the conditional expectation $E_{t}^{i,k,h}$ and we will get
\begin{equation*}
\begin{aligned}
&|\Delta Y_{t}^{i,k,h}|^{p}+\frac{1}{2}p(p-1)E_{t}^{i,k,h}\int_{t}^{T}|\Delta Y_{s}^{i,k,h}|^{p-2}|\Delta Z_{s}^{i,k,h}|^{2}ds\\
&=E_{t}^{i,k,h}\int_{t}^{T}p|\Delta Y_{s}^{i,k,h}|^{p-2}\Delta Y_{s}^{i,k,h}\alpha_{s}^{i,k,h}\Delta Y_{s}^{k,h}ds+E_{t}^{i,k,h}\int_{t}^{T}p|\Delta Y_{s}^{i,k,h}|^{p-2}\Delta Y_{s}^{i,k,h}d\Delta K_{s}^{i,k,h}.
\end{aligned}
\end{equation*}
It is easy to check that for any $p\geq2$,
\begin{align*}
&E_{t}^{i,k,h}\int_{t}^{T}|\Delta Y_{s}^{i,k,h}|^{p-2}\Delta Y_{s}^{i,k,h}d\Delta K_{s}^{i,k,h}\\
&=E_{t}^{i,k,h}\int_{t}^{T}|\Delta Y_{s}^{i,k,h}|^{p-2}(Y_{s}^{i,(k)}-S_{s}^{i})dK_{s}^{i,(k)}-E_{t}^{i,k,h}\int_{t}^{T}|\Delta Y_{s}^{i,k,h}|^{p-2}(Y_{s}^{i,(k)}-S_{s}^{i})dK_{s}^{i,(h)}\\
&\quad-E_{t}^{i,k,h}\int_{t}^{T}|\Delta Y_{s}^{i,k,h}|^{p-2}(Y_{s}^{i,(h)}-S_{s}^{i})dK_{s}^{i,(k)}+E_{t}^{i,k,h}\int_{t}^{T}|\Delta Y_{s}^{i,k,h}|^{p-2}(Y_{s}^{i,(h)}-S_{s}^{i})dK_{s}^{i,(h)}\\
&\leq-E_{t}^{i,k,h}\int_{t}^{T}|\Delta Y_{s}^{i,k,h}|^{p-2}(Y_{s}^{i,(k)}-S_{s}^{i})dK_{s}^{i,(h)}-E_{t}^{i,k,h}\int_{t}^{T}|\Delta Y_{s}^{i,k,h}|^{p-2}(Y_{s}^{i,(h)}-S_{s}^{i})dK_{s}^{i,(k)}\\
&=E_{t}^{i,k,h}\int_{t}^{T}-|(Y_{s}^{i,(k)}-S_{s}^{i})-(Y_{s}^{i,(h)}-S_{s}^{i})|^{p-2}(Y_{s}^{i,(k)}-S_{s}^{i})h(Y_{s}^{i,(h)}-S_{s}^{i})^{-}ds\\
&\quad+E_{t}^{i,k,h}\int_{t}^{T}-|(Y_{s}^{i,(h)}-S_{s}^{i})-(Y_{s}^{i,(k)}-S_{s}^{i})|^{p-2}(Y_{s}^{i,(h)}-S_{s}^{i})k(Y_{s}^{i,(k)}-S_{s}^{i})^{-}ds\\
&\leq E_{t}^{i,k,h}\int_{t}^{T}k|(Y_{s}^{i,(k)}-S_{s}^{i})^{-}|^{p-1}(Y_{s}^{i,(h)}-S_{s}^{i})^{-}ds\\
&\quad+E_{t}^{i,k,h}\int_{t}^{T}k|(Y_{s}^{i,(h)}-S_{s}^{i})^{-}|^{p-1}(Y_{s}^{i,(k)}-S_{s}^{i})^{-}ds\\
&\quad+E_{t}^{i,k,h}\int_{t}^{T}h|(Y_{s}^{i,(h)}-S_{s}^{i})^{-}|^{p-1}(Y_{s}^{i,(k)}-S_{s}^{i})^{-}ds\\
&\quad+E_{t}^{i,k,h}\int_{t}^{T}h|(Y_{s}^{i,(k)}-S_{s}^{i})^{-}|^{p-1}(Y_{s}^{i,(h)}-S_{s}^{i})^{-}ds,
\end{align*}
and
\begin{align*}
E_{t}^{i,k,h}\int_{t}^{T}|\Delta Y_{s}^{i,k,h}|^{p-2}\Delta Y_{s}^{i,k,h}\alpha_{s}^{i,k,h}\Delta Y_{s}^{k,h}ds
\leq CE_{t}^{i,k,h}\int_{t}^{T}|\Delta Y_{s}^{k,h}|^{p}ds,
\end{align*}
then it holds that
\begin{align*}
|\Delta Y_{t}^{i,k,h}|^{p}
\leq&pCE_{t}^{i,k,h}\int_{t}^{T}|\Delta Y_{s}^{k.h}|^{p}ds\\
&+pE_{t}^{i,k,h}\int_{0}^{T}(k+h)|(Y_{s}^{i,(k)}-S_{s}^{i})^{-}|^{p-1}(Y_{s}^{i,(h)}-S_{s}^{i})^{-}ds\\
&+pE_{t}^{i,k,h}\int_{0}^{T}(k+h)|(Y_{s}^{i,(h)}-S_{s}^{i})^{-}|^{p-1}(Y_{s}^{i,(k)}-S_{s}^{i})^{-}ds.
\end{align*}
We note that there exists $p(i,k,h)>1$, such that
\begin{equation*}
\begin{aligned}
E_{\nu}\left(\frac{dP^{i,k,h}}{dP}\right)^{p(i,k,h)}<\infty
\end{aligned}
\end{equation*}
from \cite[Theorem 3.1]{20}, where $E_{\nu}$ is a conditional expectation for an arbitrary stopping time $\nu$. Let $q(i,k,h)>1$, such that
\begin{equation*}
\begin{aligned}
\frac{1}{p(i,k,h)}+\frac{1}{q(i,k,h)}=1.
\end{aligned}
\end{equation*}
By Doob's $L^{p}$-inequality, H\"older's inequality and Young's inequality, we show that for any $p\geq2,~k,h\geq0,~i=1,\cdots,n$,
\begin{align*}
&E^{i,k,h}\left[\sup_{t\in[0,T]}E_{t}^{i,k,h}\int_{t}^{T}|\Delta Y_{s}^{k,h}|^{p}ds\right]\\
&\leq 2E^{i,k,h}\left[\left(\int_{0}^{T}|\Delta Y_{s}^{k,h}|^{p}ds\right)^{2}\right]^{\frac{1}{2}}\\
&=2E\left[\left(\frac{dP^{i,k,h}}{dP}\right)\left(\int_{0}^{T}|\Delta Y_{s}^{k,h}|^{p}ds\right)^{2}\right]^{\frac{1}{2}}\\
&\leq2E\left[\left(\frac{dP^{i,k,h}}{dP}\right)^{p(i,k,h)}\right]^{\frac{1}{2p(i,k,h)}}E\left[\left(\int_{0}^{T}|\Delta Y_{s}^{k,h}|^{p}ds\right)^{2q(i,k,h)}\right]^{\frac{1}{2q(i,k,h)}}\\
&\leq2E\left[\left(\frac{dP^{i,k,h}}{dP}\right)^{p(i,k,h)}\right]^{\frac{1}{2p(i,k,h)}}E\left[\left(\int_{0}^{T}|\Delta Y_{s}^{k,h}|^{2pq(i,k,h)}ds\right)\right]^{\frac{1}{2q(i,k,h)}},
\end{align*}
and
\begin{align*}
&E^{i,k,h}\left[\sup_{t\in[0,T]}E_{t}^{i,k,h}\int_{t}^{T}h|(Y_{s}^{i,(k)}-S_{s}^{i})^{-}|^{p-1}(Y_{s}^{i,(h)}-S_{s}^{i})^{-}ds\right]\\
&= E^{i,k,h}\left[\sup_{t\in[0,T]}E_{t}^{i,k,h}\int_{t}^{T}|(Y_{s}^{i,(k)}-S_{s}^{i})^{-}|^{p-1}dK_{s}^{i,(h)}\right]\\
&\leq 2E^{i,k,h}\left[\left(\int_{0}^{T}|(Y_{s}^{i,(k)}-S_{s}^{i})^{-}|^{p-1}dK_{s}^{i,(h)}\right)^{2}\right]^{\frac{1}{2}}\\
&\leq 2E^{i,k,h}\left[\sup_{t\in[0,T]}|(Y_{t}^{i,(k)}-S_{t}^{i})^{-}|^{2(p-1)}(K_{T}^{i,(h)})^{2}\right]^{\frac{1}{2}}\\
&= 2E\left[\left(\frac{dP^{i,k,h}}{dP}\right)\sup_{t\in[0,T]}|(Y_{t}^{i,(k)}-S_{t}^{i})^{-}|^{2(p-1)}(K_{T}^{i,(h)})^{2}\right]^{\frac{1}{2}}\\
&\leq 2B_{2q(i,k,h)}E\left[\left(\frac{dP^{i,k,h}}{dP}\right)^{p(i,k,h)}\right]^{\frac{1}{2p(i,k,h)}}E\left[\sup_{t\in[0,T]}|(Y_{t}^{i,(k)}-S_{t}^{i})^{-}|^{2(p-1)q(i,k,h)}\right]^{\frac{1}{2q(i,k,h)}},
\end{align*}
and
\begin{align*}
&E^{i,k,h}\left[\sup_{t\in[0,T]}E_{t}^{i,k,h}\int_{t}^{T}k|(Y_{s}^{i,(k)}-S_{s}^{i})^{-}|^{p-1}(Y_{s}^{i,(h)}-S_{s}^{i})^{-}ds\right]\\
&= E^{i,k,h}\left[\sup_{t\in[0,T]}E_{t}^{i,k,h}\int_{t}^{T}|(Y_{s}^{i,(k)}-S_{s}^{i})^{-}|^{p-2}|(Y_{s}^{i,(h)}-S_{s}^{i})^{-}|dK_{s}^{i,(k)}\right]\\
&\leq \frac{p-2}{p-1}E^{i,k,h}\left[\sup_{t\in[0,T]}E_{t}^{i,k,h}\int_{t}^{T}|(Y_{s}^{i,(k)}-S_{s}^{i})^{-}|^{p-1}dK_{s}^{i,(k)}\right]\\
&\quad+\frac{1}{p-1}E^{i,k,h}\left[\sup_{t\in[0,T]}E_{t}^{i,k,h}\int_{t}^{T}|(Y_{s}^{i,(h)}-S_{s}^{i})^{-}|^{p-1}dK_{s}^{i,(k)}\right]\\
&\leq \frac{2(p-2)B_{2q(i,k,h)}}{p-1}E\left[\left(\frac{dP^{i,k,h}}{dP}\right)^{p(i,k,h)}\right]^{\frac{1}{2p(i,k,h)}}E\left[\sup_{t\in[0,T]}|(Y_{t}^{i,(k)}-S_{t}^{i})^{-}|^{2(p-1)q(i,k,h)}\right]^{\frac{1}{2q(i,k,h)}}\\
&\quad+\frac{2B_{2q(i,k,h)}}{p-1}E\left[\left(\frac{dP^{i,k,h}}{dP}\right)^{p(i,k,h)}\right]^{\frac{1}{2p(i,k,h)}}E\left[\sup_{t\in[0,T]}|(Y_{t}^{i,(h)}-S_{t}^{i})^{-}|^{2(p-1)q(i,k,h)}\right]^{\frac{1}{2q(i,k,h)}},
\end{align*}
which implies
\begin{equation*}
\begin{aligned}
E^{i,k,h}&\left[\sup_{t\in[0,T]}|\Delta Y_{t}^{i,k,h}|^{p}\right]
\leq 2pCE\left[\left(\frac{dP^{i,k,h}}{dP}\right)^{p(i,k,h)}\right]^{\frac{1}{2p(i,k,h)}}E\left[\left(\int_{0}^{T}|\Delta Y_{s}^{k,h}|^{2pq(i,k,h)}ds\right)\right]^{\frac{1}{2q(i,k,h)}}\\
&+4pB_{2q(i,k,h)}E\left[\left(\frac{dP^{i,k,h}}{dP}\right)^{p(i,k,h)}\right]^{\frac{1}{2p(i,k,h)}}E\left[\sup_{t\in[0,T]}|(Y_{t}^{i,(k)}-S_{t}^{i})^{-}|^{2(p-1)q(i,k,h)}\right]^{\frac{1}{2q(i,k,h)}}\\
&+4pB_{2q(i,k,h)}E\left[\left(\frac{dP^{i,k,h}}{dP}\right)^{p(i,k,h)}\right]^{\frac{1}{2p(i,k,h)}}E\left[\sup_{t\in[0,T]}|(Y_{t}^{i,(h)}-S_{t}^{i})^{-}|^{2(p-1)q(i,k,h)}\right]^{\frac{1}{2q(i,k,h)}}.
\end{aligned}
\end{equation*}
We note that there exists $p'(i,k,h)>1$, such that
\begin{equation*}
\begin{aligned}
E^{i,k,h}_{\nu}\left(\frac{dP}{dP^{i,k,h}}\right)^{p'(i,k,h)}<\infty
\end{aligned}
\end{equation*}
from \cite[Theorem 3.1]{20}, where $E^{i,k,h}_{\nu}$ is a conditional expectation for an arbitrary stopping time $\nu$. Let $q'(i,k,h)>1$, such that
\begin{equation*}
\begin{aligned}
\frac{1}{p'(i,k,h)}+\frac{1}{q'(i,k,h)}=1.
\end{aligned}
\end{equation*}
For $p\geq2$,
\begin{align*}
E\left[\sup_{t\in[0,T]}|\Delta Y_{t}^{i,k,h}|^{p}\right]&=E^{i,k,h}\left[\left(\frac{dP}{dP^{i,k,h}}\right)\sup_{t\in[0,T]}|\Delta Y_{t}^{i,k,h}|^{p}\right]\\
&\leq E^{i,k,h}\left[\left(\frac{dP}{dP^{i,k,h}}\right)^{p'(i,k,h)}\right]^{\frac{1}{p'(i,k,h)}}E^{i,k,h}\left[\sup_{t\in[0,T]}|\Delta Y_{t}^{i,k,h}|^{pq'(i,k,h)}\right]^{\frac{1}{q'(i,k,h)}}
\end{align*}
which means $Y^{i,(k)}$ is convergent in $\mathcal{S}^{p}(\mathbb{R})$ for any $p\geq2$. So $\|Y^{(k)}-Y\|_{\mathcal{S}^{p}(\mathbb{R}^{n})}\rightarrow0$ as $k\rightarrow\infty$ for any $p\geq2$.\\
Use It\^o's formula to $|\Delta Y_{t}^{i,k,h}|^{2}$ and we will get
\begin{align*}
&|\Delta Y_{t}^{i,k,h}|^{2}+\int_{t}^{T}|\Delta Z_{s}^{i,k,h}|^{2}ds\\
&=2\int_{t}^{T}\Delta Y_{s}^{i,k,h}(\alpha_{s}^{i,k,h}\Delta Y_{s}^{k,h}+\beta_{s}^{i,k,h}\Delta Z_{s}^{i,k,h})ds+2\int_{t}^{T}\Delta Y_{s}^{i,k,h}d\Delta K_{s}^{i,k,h}\\
&\quad-2\int_{t}^{T}\Delta Y_{s}^{i,k,h}\Delta Z_{s}^{i,k,h}dW_{s}.
\end{align*}
Let t=0, we get
\begin{align*}
\int_{0}^{T}|\Delta Z_{s}^{i,k,h}|^{2}ds
\leq&2\int_{0}^{T}\Delta Y_{s}^{i,k,h}(\alpha_{s}^{i,k,h}\Delta Y_{s}^{k,h}+\beta_{s}^{i,k,h}\Delta Z_{s}^{i,k,h})ds\\
&+2\int_{0}^{T}\Delta Y_{s}^{i,k,h}d\Delta K_{s}^{i,k,h}-2\int_{0}^{T}\Delta Y_{s}^{i,k,h}\Delta Z_{s}^{i,k,h}dW_{s}\\
\leq&2CT\sup_{t\in[0,T]}|\Delta Y_{t}^{k,h}|^{2}+6C\sup_{t\in[0,T]}|\Delta Y_{t}^{k,h}|\int_{0}^{T}(1+|Z_{s}^{i,(k)}|^{2}+|Z_{s}^{i,(h)}|^{2})ds\\
&+2\sup_{t\in[0,T]}|\Delta Y_{t}^{k,h}|(|K_{T}^{(k)}|+|K_{T}^{(h)}|)-2\int_{0}^{T}\Delta Y_{s}^{i,k,h}\Delta Z_{s}^{i,k,h}dW_{s}.
\end{align*}
With $C_{r}$ inequality and H\"older's inequality, for any $p\geq2$,
\begin{align*}
&E[(\int_{0}^{T}|\Delta Z_{s}^{i,k,h}|^{2}ds)^{\frac{p}{2}}]\\
&\leq C_{\frac{p}{2}}\left\{(2CT)^{\frac{p}{2}}E\left[\sup_{t\in[0,T]}|\Delta Y_{t}^{k,h}|^{p}\right]\right.\\
&\quad+(6C)^{\frac{p}{2}}E\left[\sup_{t\in[0,T]}|\Delta Y_{t}^{k,h}|^{\frac{p}{2}}\left(\int_{0}^{T}(1+|Z_{s}^{i,(k)}|^{2}+|Z_{s}^{i,(h)}|^{2})ds\right)^{\frac{p}{2}}\right]\\
&\quad+\left.2^{\frac{p}{2}}E\left[\sup_{t\in[0,T]}|\Delta Y_{t}^{k,h}|^{\frac{p}{2}}(|K_{T}^{(k)}|^{\frac{p}{2}}+|K_{T}^{(h)}|^{\frac{p}{2}})\right]+2^{\frac{p}{2}}E\left[\left(\int_{0}^{T}|\Delta Y_{s}^{i,k,h}\Delta Z_{s}^{i,k,h}|^{2}ds\right)^{\frac{p}{4}}\right]\right\}\\
&\leq C'_{\frac{p}{2}}\left\{E\left[\sup_{t\in[0,T]}|\Delta Y_{t}^{k,h}|^{p}\right]+E\left[\sup_{t\in[0,T]}|\Delta Y_{t}^{k,h}|^{p}\right]^{\frac{1}{2}}E\left[\left(\int_{0}^{T}(1+|Z_{s}^{i,(n)}|^{2}+|Z_{s}^{i,(m)}|^{2})ds\right)^{p}\right]^{\frac{1}{2}}\right.\\
&\quad+E\left[\sup_{t\in[0,T]}|\Delta Y_{t}^{k,h}|^{p}\right]^{\frac{1}{2}}(E[|K_{T}^{(k)}|^{p}]^{\frac{1}{2}}+E[|K_{T}^{(h)}|^{p}]^{\frac{1}{2}})\\
&\quad\left.+E\left[\sup_{t\in[0,T]}|\Delta Y_{t}^{k,h}|^{p}\right]^{\frac{1}{2}}E\left[\left(\int_{0}^{T}|Z_{s}^{i,(k)}|^{2}+|Z_{s}^{i,(h)}|^{2}ds\right)^{\frac{p}{2}}\right]^{\frac{1}{2}}\right\}.
\end{align*}
Since for any $p\geq2,k\geq0$,
\begin{equation*}
\begin{aligned}
\lim\limits_{k,h\rightarrow\infty}E\left[\sup\limits_{t\in[0,T]}|\Delta Y_{t}^{k,h}|^{p}\right]=0,~ E\left(\int_{0}^{T}|Z_{s}^{i,(k)}|^{2}ds\right)^{\frac{p}{2}}<\infty,~E|K_{T}^{(k)}|^{p}<\infty,
\end{aligned}
\end{equation*}
then $Z^{i,(k)}$ is convergent in $\mathcal{H}^{p}(\mathbb{R}^{d})$. So there exists a process $Z\in\mathcal{H}^{p}(\mathbb{R}^{n\times d})$, such that $\|Z^{(k)}-Z\|_{\mathcal{H}^{p}(\mathbb{R}^{n\times d})}\rightarrow0$ as $k\rightarrow\infty$ for any $p\geq2$.\\
The definition of $K_{t}^{(k)}$ shows it is a nondecreasing process and
\begin{align*}
\Delta K_{t}^{i,k,h}=\Delta Y_{0}^{i,k,h}-\Delta Y_{t}^{i,k,h}-\int_{0}^{t}(\alpha_{s}^{i,k,h}\Delta Y_{s}^{k,h}+\beta_{s}^{i,k,h}\Delta Z_{s}^{i,k,h})ds+\int_{0}^{t}\Delta Z_{s}^{i,k,h}dW_{s}.
\end{align*}
For any $p\geq2$,
\begin{align*}
E\left[\sup_{t\in[0,T]}|\Delta K_{t}^{i,k,h}|^{p}\right]
\leq&C_{p}\left\{2E\left[\sup_{t\in[0,T]}|\Delta Y_{t}^{i,k,h}|^{p}\right]
+(CT)^{p}E\left[\sup_{t\in[0,T]}|\Delta Y_{t}^{i,k,h}|^{p}\right]\right.\\
&+C^{p}E\left[\left(\int_{0}^{T}(1+|Z_{s}^{i,(k)}|+|Z_{s}^{i,(h)}|)|\Delta Z_{s}^{i,k,h}|ds\right)^{p}\right]\\
&+\left.E\left[\left(\int_{0}^{T}|\Delta Z_{s}^{i,k,h}|^{2}ds\right)^{\frac{p}{2}}\right]\right\}\\
\leq&C'_{p}\left\{E\left[\sup_{t\in[0,T]}|\Delta Y_{t}^{i,k,h}|^{p}\right]+E\left[\left(\int_{0}^{T}|\Delta Z_{s}^{i,k,h}|^{2}ds\right)^{\frac{p}{2}}\right]\right.\\
&+\left.E\left[\left(\int_{0}^{T}(1+|Z_{s}^{i,(k)}|+|Z_{s}^{i,(h)}|)^{2}ds\right)^{p}\right]^{\frac{1}{2}}E\left[\left(\int_{0}^{T}|\Delta Z_{s}^{i,k,h}|^{2}ds\right)^{p}\right]^{\frac{1}{2}}\right\},
\end{align*}
which implies $K^{i,(k)}$ is convergent in $\mathcal{S}^{p}(\mathbb{R}^{n})$. So there exists a continuous nondecreasing process $K\in\mathcal{S}^{p}(\mathbb{R}^{n})$, such that $\|K^{(k)}-K\|_{\mathcal{S}^{p}(\mathbb{R}^{n})}\rightarrow0$ as $k\rightarrow\infty$ for any $p\geq2$,
and
\begin{align*}
K_{0}=0,~\int_{0}^{T}(Y_{t}^{i}-S_{t}^{i})dK_{t}^{i}=\lim\limits_{k\rightarrow\infty}\int_{0}^{T}k(Y_{t}^{i,(k)}-S_{t}^{i})(Y_{t}^{i,(k)}-S_{t}^{i})^{-}ds=0.
\end{align*}
Therefore let $k\rightarrow\infty$ in (2), we have
\begin{align*}
Y_{t}^{i}=\xi^{i}+\int_{t}^{T}f^{i}(s,Y_{s},Z_{s}^{i})ds+K_{T}^{i}-K_{t}^{i}-\int_{t}^{T}Z_{s}^{i}dW_{s}.
\end{align*}
Since $\|Y^{(k)}\|_{\mathcal{S}^{\infty}(\mathbb{R}^{n})}\leq C^{*}$ and $\|Y^{(k)}-Y\|_{\mathcal{S}^{p}(\mathbb{R}^{n})}\rightarrow0$ as $k\rightarrow\infty$ for any $p\geq2$, we have $\|Y\|_{\mathcal{S}^{\infty}(\mathbb{R}^{n})}\leq C^{*}$. Then similar to the inequality (3), we have
\begin{align*}
4C^{2}e^{-4CC^{*}}E_{t}\int_{t}^{T}|Z_{s}^{i}|^{2}ds
\leq e^{4C^{2}}+4C^{2}(1+C^{*})Te^{4CC^{*}},
\end{align*}
which implies $Z\in BMO$.\\
\\
$\mathbf{Step\ 4}$: We show the uniqueness of the solution.\\
Let $(Y,Z,K)$ and $(\overline{Y},\overline{Z},\overline{K})$ be two solutions of the RBSDE (1).\\
Set
\begin{align*}
&\widehat{Y}=Y-\overline{Y},~\widehat{Y}^{i}=Y^{i}-\overline{Y}^{i},~\widehat{Z}^{i}=Z^{i}-\overline{Z}^{i},~\widehat{K}^{i}=K^{i}-\overline{K}^{i},\\
&f^{i}(t,Y_{t},Z_{t}^{i})-f^{i}(t,\overline{Y}_{t},\overline{Z}_{t}^{i})=\alpha^{i}_{t} \widehat{Y}_{t}+\beta^{i}_{t}\widehat{Z}_{t}^{i},\\
&|\alpha^{i}_{t}|\leq C,~|\beta^{i}_{t}|\leq C(1+|Z_{t}^{i}|+|\overline{Z}_{t}^{i}|).
\end{align*}
Then $W_{t}^{i}:=W_{t}-\int_{0}^{t}\beta_{s}^{i}ds$ is a Brownian motion under the equivalent probability measure $P^{i}$ defined by
\begin{align*}
dP^{i}:=\mathcal{E}(\beta^{i}\cdot W)_{0}^{T}dP.
\end{align*}
Use It\^o's formula to $|\widehat{Y}_{t}^{i}|^{2}$ under the conditional expectation $E_{t}^{i}$ and we will get
\begin{align*}
|\widehat{Y}_{t}^{i}|^{2}+E_{t}^{i}\int_{t}^{T}|\widehat{Z}_{s}^{i}|^{2}ds=2E_{t}^{i}\int_{t}^{T}\widehat{Y}_{s}^{i}\alpha_{s}^{i}\widehat{Y}_{s}ds+2E_{t}^{i}\int_{t}^{T}\widehat{Y}_{t}^{i}d\widehat{K}_{s}^{i}.
\end{align*}
Since
\begin{align*}
2E_{t}^{i}\int_{t}^{T}\widehat{Y}_{s}^{i}\alpha_{s}^{i}\widehat{Y}_{s}ds
\leq CE_{t}^{i}\int_{t}^{T}2|\widehat{Y}_{s}^{i}||\widehat{Y}_{s}|ds
\leq CE_{t}^{i}\int_{t}^{T}|\widehat{Y}_{s}^{i}|^{2}ds+C\sum\limits_{j=1}^{n}E_{t}^{i}\int_{t}^{T}|\widehat{Y}_{s}^{j}|^{2}ds
\end{align*}
and
\begin{align*}
E_{t}^{i}\int_{t}^{T}\widehat{Y}_{s}^{i}d\widehat{K}_{s}^{i}=&E_{t}^{i}\int_{t}^{T}(Y_{s}^{i}-S_{s}^{i})dK_{s}^{i}-E_{t}^{i}\int_{t}^{T}(Y_{s}^{i}-S_{s}^{i})d\overline{K}_{s}^{i}\\
&-E_{t}^{i}\int_{t}^{T}(\overline{Y}_{s}^{i}-S_{s}^{i})dK_{s}^{i}+E_{t}^{i}\int_{t}^{T}(\overline{Y}_{s}^{i}-S_{s}^{i})d\overline{K}_{s}^{i}\\
=&-E_{t}^{i}\int_{t}^{T}(Y_{s}^{i}-S_{s}^{i})d\overline{K}_{s}^{i}-E_{t}^{i}\int_{t}^{T}(\overline{Y}_{s}^{i}-S_{s}^{i})dK_{s}^{i}\\
\leq&0,
\end{align*}
we have
\begin{align*}
&\|\sup_{t\leq s\leq T}|\widehat{Y}_{s}^{i}|^{2}\|_{\infty}\leq C(T-t)(\|\sup_{t\leq s\leq T}|\widehat{Y}_{s}^{i}|^{2}\|_{\infty}+\sum\limits_{j=1}^{n}\|\sup_{t\leq s\leq T}|\widehat{Y}_{s}^{j}|^{2}\|_{\infty}).
\end{align*}
Sum up the equation above respect to $i$, we can obtain
\begin{align*}
\sum\limits_{i=1}^{n}\|\sup_{t\leq s\leq T}|\widehat{Y}_{s}^{i}|^{2}\|_{\infty}\leq C(T-t)(\sum\limits_{i=1}^{n}\|\sup_{t\leq s\leq T}|\widehat{Y}_{s}^{i}|^{2}\|_{\infty}+n\sum\limits_{j=1}^{n}\|\sup_{t\leq s\leq T}|\widehat{Y}_{s}^{j}|^{2}\|_{\infty}),
\end{align*}
which is equivalent to
\begin{align*}
[1-C(1+n)(T-t)]\sum\limits_{i=1}^{n}\|\sup_{t\leq s\leq T}|\widehat{Y}_{s}^{i}|^{2}\|_{\infty}\leq 0.
\end{align*}
Let $\epsilon=1-C(1+n)(T-t)$. When $T-t=\varepsilon>0$ is small enough, we have $0<\epsilon<1$ and
\begin{align*}
\sum\limits_{i=1}^{n}\|\sup\limits_{t\leq s\leq T}|\widehat{Y}_{s}^{i}|^{2}\|_{\infty}\leq 0.
\end{align*}
Since the terminal time $T$ is finite, after a finite step iteration, we have
\begin{align*}
\sum\limits_{i=1}^{n}\|\sup\limits_{0\leq s\leq T}|\widehat{Y}_{s}^{i}|^{2}\|_{\infty}\leq 0,
\end{align*}
which means $\|\widehat{Y}\|_{\mathcal{S}^{\infty}(\mathbb{R}^{n})}=0$ and then $Y_{t}=\overline{Y}_{t},~t\in[0,T]$.\\
Similar to the result in the step 3, for any $p\geq2$, we have
\begin{align*}
E[(\int_{0}^{T}|\widehat{Z}_{s}^{i}|^{2}ds)^{\frac{p}{2}}]\leq& C'_{\frac{p}{2}}\left\{E\left[\sup_{t\in[0,T]}|\widehat{Y}_{t}|^{p}\right]+E\left[\sup_{t\in[0,T]}|\widehat{Y}_{t}|^{p}\right]^{\frac{1}{2}}E\left[\left(\int_{0}^{T}(1+|Z_{s}^{i}|^{2}+|\overline{Z}_{s}^{i}|^{2})ds\right)^{p}\right]^{\frac{1}{2}}\right.\\
&+E\left[\sup_{t\in[0,T]}|\widehat{Y}_{t}|^{p}\right]^{\frac{1}{2}}(E[|K_{T}|^{p}]^{\frac{1}{2}}+E[|\overline{K}_{T}|^{p}]^{\frac{1}{2}})\\
&+\left.E\left[\sup_{t\in[0,T]}|\widehat{Y}_{t}|^{p}\right]^{\frac{1}{2}}E\left[\left(\int_{0}^{T}|Z_{s}^{i}|^{2}+|\overline{Z}_{s}^{i}|^{2}ds\right)^{\frac{p}{2}}\right]^{\frac{1}{2}}\right\}\\
=&0.
\end{align*}
At last, for any $p\geq2$, we have
\begin{align*}
E\left[\sup_{t\in[0,T]}|\widehat{K}_{t}^{i}|^{p}\right]
\leq&C'_{p}\{E\left[\sup_{t\in[0,T]}|\widehat{Y}_{t}^{i}|^{p}\right]+E\left[\left(\int_{0}^{T}|\widehat{Z}_{s}^{i}|^{2}ds\right)^{\frac{p}{2}}\right]\}\\
&+E\left[\left(\int_{0}^{T}(1+|Z_{s}^{i}|+|\overline{Z}_{s}^{i}|)^{2}ds\right)^{p}\right]^{\frac{1}{2}}E\left[\left(\int_{0}^{T}|\widehat{Z}_{s}^{i}|^{2}ds\right)^{p}\right]^{\frac{1}{2}}\}\\
=&0.
\end{align*}
We finish the uniqueness part.\\
\\
In conclusion, the RBSDE (1) has a unique solution $(Y,Z,K)$ such that $(Y,Z\cdot W,K)\in\mathcal{S}^{\infty}(\mathbb{R}^{n})\times BMO\times \mathcal{S}^{p}(\mathbb{R}^{n})$ for any $ p\geq2$.
\end{proof}

\begin{theorem}
Let $(Y,Z,K)$ be a solution of the RBSDE (1) satisfying $(\mathscr{A}1)-(\mathscr{A}4)$. Then for $\forall t\in[0,T],~i=1,\cdots,n$,
\begin{align*}
Y_{t}^{i}=\esssup\limits_{\nu\in\mathcal{T}_{t}}E_{t}[\int_{t}^{\nu}f^{i}(s,Y_{s},Z_{s}^{i})ds+S_{\nu}^{i}1_{\{\nu<T\}}+\xi^{i}1_{\{\nu=T\}}]
\end{align*}
where $\mathcal{T}$ is the set of all stopping times dominated by $T$ and $\mathcal{T}_{t}=\{\nu\in\Gamma;t\leq\nu\leq T\}$.
\end{theorem}
\begin{proof}
Let $\nu\in\mathcal{T}_{t}$ and take the conditional expectation, we have
\begin{align*}
Y_{t}^{i}=&E_{t}[\int_{t}^{\nu}f^{i}(s,Y_{s},Z_{s}^{i})ds+Y_{\nu}^{i}+K_{\nu}^{i}-K_{t}^{i}]\\
\geq&E_{t}[\int_{t}^{\nu}f^{i}(s,Y_{s},Z_{s}^{i})ds+S_{\nu}^{i}1_{\{\nu<T\}}+\xi^{i}1_{\{\nu=T\}}].
\end{align*}
Then we choose a suitable stopping time to get the equation. Let $D_{t}=\inf\{t\leq s\leq T;Y_{s}^{i}=S_{s}^{i}\}$, then
\begin{align*}
Y_{s}^{i}>S_{s}^{i},~t\leq s<D_{t}
\end{align*}
and
\begin{align*}
K_{D_{t}}^{i}=K_{t}^{i}.
\end{align*}
It follows that
\begin{align*}
Y_{t}^{i}=E_{t}[\int_{t}^{D_{t}}f^{i}(s,Y_{s},Z_{s}^{i})ds+S_{D_{t}}^{i}1_{\{D_{t}<T\}}+\xi^{i}1_{\{D_{t}=T\}}].
\end{align*}
Hence, the results follow.
\end{proof}

\begin{theorem}
Let $(Y,Z,K)$ (resp. $(\overline{Y},\overline{Z},\overline{K})$) be the solution of RBSDEs (1) associated to  $(\xi,f,S)$ (resp. $(\overline{\xi},\overline{f},\overline{S})$). Assume that the terminal conditions, generators and barriers $(\xi,f,S)$ and $(\overline{\xi},\overline{f},\overline{S})$ satisfying $(\mathscr{A}1)-(\mathscr{A}4)$. If it holds that $\xi\leq\overline{\xi},~S_{t}\leq\overline{S}_{t}$ and for any $t\in[0,T],~i=1,\cdots,n$,
\begin{align*}
f^{i}(t,y,z^{i})\leq \overline{f}^{i}(t,\overline{y},z^{i})
\end{align*}
for any $z^{i}\in\mathbb{R}^{d}$ and $y,\overline{y}\in\mathbb{R}^{n}$ satisfying $y^{i}=\overline{y}^{i},~y^{j}\leq\overline{y}^{j},j\neq i$, we have
\begin{align*}
Y_{t}\leq \overline{Y}_{t},~t\in[0,T].
\end{align*}
\end{theorem}
\begin{proof}\\
For each $k\in\mathbb{N}$, we consider two BSDEs:
\begin{align*}
&Y_{t}^{i,(k)}=\xi^{i}+\int_{t}^{T}f^{i}(s,Y_{s}^{(k)},Z_{s}^{i,(k)})ds+k\int_{t}^{T}(Y_{s}^{i,(k)}-S_{s}^{i})^{-}ds-\int_{t}^{T}Z_{s}^{i,(k)}dW_{s},\\
&\overline{Y}_{t}^{i,(k)}=\overline{\xi}^{i}+\int_{t}^{T}\overline{f}^{i}(s,\overline{Y}_{s}^{(k)},\overline{Z}_{s}^{i,(k)})ds+k\int_{t}^{T}(\overline{Y}_{s}^{i,(k)}-\overline{S}_{s}^{i})^{-}ds-\int_{t}^{T}\overline{Z}_{s}^{i,(k)}dW_{s}.
\end{align*}
It follows from Theorem 2.1 that each equation above has a solution and $(Y^{(k)},Z^{(k)}\cdot W),(\overline{Y}^{(k)},~\overline{Z}^{(k)}\cdot W)\in\mathcal{S}^{\infty}(\mathbb{R}^{n})\times BMO$. We define
\begin{align*}
&f_{k}^{i}(t,y,z^{i})=f^{i}(t,y,z^{i})+k(y^{i}-S_{s}^{i})^{-},\\
&\overline{f}_{k}^{i}(t,\overline{y},z^{i})=f^{i}(t,\overline{y},z^{i})+k(\overline{y}^{i}-\overline{S}_{s}^{i})^{-}.
\end{align*}
Since $S_{t}^{i}\leq\overline{S}_{t}^{i}$, we have $(y^{i}-S_{t}^{i})^{-}\leq(\overline{y}^{i}-\overline{S}_{t}^{i})^{-}$ and then
\begin{align*}
f_{k}^{i}(t,y,z^{i})\leq \overline{f}_{k}^{i}(t,\overline{y},z^{i}),~t\in[0,T]
\end{align*}
for any $z^{i}\in\mathbb{R}^{d}$ and $y,\overline{y}\in\mathbb{R}^{n}$ satisfying $y^{i}=\overline{y}^{i},y^{j}\leq\overline{y}^{j},j\neq i$.\\
Using \cite[Theorem 2.2]{19}, it holds that
\begin{align*}
Y_{t}^{(k)}\leq \overline{Y}_{t}^{(k)},~t\in[0,T],~k\geq 0.
\end{align*}
Afterwards we know that
\begin{align*}
Y_{t}=\lim_{k\rightarrow\infty}Y_{t}^{(k)}\leq \lim_{k\rightarrow\infty}\overline{Y}_{t}^{(k)}=\overline{Y}_{t},~t\in[0,T].
\end{align*}
Hence, we obtain the results.
\end{proof}

\end{document}